\documentclass[12pt,a4paper]{amsart}
\usepackage[utf8]{inputenc}
\usepackage{amsmath, amsthm, amssymb, amsfonts, mathabx}
\usepackage{geometry}
\usepackage[rightcaption]{sidecap}
\usepackage{enumitem}
\usepackage{hyperref}
\usepackage{stackengine} 

\newcommand*{\Z}{\mathbb{Z}}%
\usepackage{tkz-graph}
\usepackage{tkz-berge}

\theoremstyle{plain}

\newtheorem{proposition}{Proposition}

\theoremstyle{definition}
\newtheorem*{acknowledgement}{Acknowledgement}

\begin{document}
\bibliographystyle{plain}

\title{Discrete logarithm problem in some families of sandpile groups}
\author{K. Dsupin}
\email{krisztian.dsupin@gmail.com}
\author{Sz. Tengely}
\email{tengely@science.unideb.hu}
\address{Institute of Mathematics, University of Debrecen \\
P. O. Box 400, 4002 Debrecen, Hungary}


\begin{abstract}
Biggs  proposed the sandpile group of certain modified wheel graphs  for cryptosystems relying on the difficulty of the discrete logarithm problem. Blackburn and independently Shokrieh showed that the discrete logarithm problem is  efficiently solvable. We study Shokrieh's method in cases of graphs such that the sandpile group is not cyclic, namely the square cycle graphs and the wheel graphs. Knowing generators of the group or the form of the pseudoinverse of the Laplacian matrix makes the problem more vulnerable. We also consider the discrete logarithm problem in case of the so-called subdivided banana graphs. In certain cases the sandpile group is cyclic and a generator is known and one can solve the discrete logarithm problem without computing the pseudoinverse of the Laplacian matrix.
\end{abstract}


\maketitle

\section{Introduction}
Finding a usable representation of a group for discrete logarithm problem (DLP) is a popular research area in cryptography. 
There are families of groups in which cases the DLP is solvable, hence should not be used to design cryptographic primitives. Ariffin and Abu \cite{AAbeta} proposed a scheme based on the infinite additive group $\mathbb{R}/\mathbb{Z}.$ Blackburn \cite{Blackburn_AAbeta} studied the security of the public key scheme and showed that it can be attacked efficiently by means of the continued fraction algorithm. 
Ili\'c and Magliveras \cite{DLP_NA} considered the DLP in non-abelian groups, that is based on a particular representation of the and a choice of generators. They showed that in $PSL(2,p)$ in certain cases the problem is easy to solve, therefore such representation of $PSL(2, p)$ and generators should not be used to design cryptographic primitives.
Sandpile groups of certain class of graphs were suggested for cryptosystems 
by Biggs in \cite{bigg_cryto}. His proposal was a modification of the wheel graph, which has a special property to form a cyclic sandpile group. Blackburn in \cite{Blackburn} and Shokrieh in \cite{monodromy_pairing} independently showed that the DLP in that case is efficiently solvable.
Moreover, Shokrieh suggested that with the modification of his method the DLP can be also solved in non-cyclic sandpile groups. The cardinality of the sandpile group is equal to the number of spanning trees of the graphs. The literature is rich in the direction of determining the structure of the sandpile group of certain families of graphs. Some examples out of these works are: wheel graphs - Biggs \cite{GWheel}; Paley graphs - Chandler, Sin and Xiang \cite{GPaley}; polygon flower graph - Chen and Mohar \cite{GPolygon}; $P_4\times C_n$ graph - Chen and Hou \cite{GP4Cn}; Möbius ladder graph - Chen, Hou and Woo \cite{GMobius}; rook's graph - Ducey, Gerhard and Watson \cite{GRook}. 

To be able to apply Shokrieh's method one needs to know generators of the given sandpile group and the pseudoinverse of the Laplacian matrix. We remark that in some classes of graphs there are results providing generators (like in case of wheel graphs \cite{GWheel}) and in case of circulant Laplacians  Cline, Plemmons and Worm \cite{generalized_inverse} give an efficient computational form. In case of general graph Laplacians there is a GPU-based method to compute the pseudoinverse by Saurabh, Varbanescu and Ranjan \cite{GPU_L}, they provided numerical examples related to graphs having about 8000 vertices.

In this article we study the method of Shokrieh in case of three families of graphs, the square cycle graphs, the wheel graphs and subdivided banana graphs.

\section{DLP in square cycle graphs}
Hou, Woo and Chen studied the $C_n^2$ graph in \cite{cn2_sandpile}. 
Let the vertex set of the cycle $C_n$ be $V = \{1, 2,\ldots, n\}.$ The square cycle graph $C_n^2$ is a graph with vertex set $V = \{1, 2,\ldots , n\},$ and each vertex $i$ is adjacent to vertices $i \pm 1, i \pm 2 \bmod{n}.$
\begin{figure}
\begin{tikzpicture}
\definecolor{cv0}{rgb}{0.0,0.0,0.0}
\definecolor{cfv0}{rgb}{1.0,1.0,1.0}
\definecolor{clv0}{rgb}{0.0,0.0,0.0}
\definecolor{cv1}{rgb}{0.0,0.0,0.0}
\definecolor{cfv1}{rgb}{1.0,1.0,1.0}
\definecolor{clv1}{rgb}{0.0,0.0,0.0}
\definecolor{cv2}{rgb}{0.0,0.0,0.0}
\definecolor{cfv2}{rgb}{1.0,1.0,1.0}
\definecolor{clv2}{rgb}{0.0,0.0,0.0}
\definecolor{cv3}{rgb}{0.0,0.0,0.0}
\definecolor{cfv3}{rgb}{1.0,1.0,1.0}
\definecolor{clv3}{rgb}{0.0,0.0,0.0}
\definecolor{cv4}{rgb}{0.0,0.0,0.0}
\definecolor{cfv4}{rgb}{1.0,1.0,1.0}
\definecolor{clv4}{rgb}{0.0,0.0,0.0}
\definecolor{cv5}{rgb}{0.0,0.0,0.0}
\definecolor{cfv5}{rgb}{1.0,1.0,1.0}
\definecolor{clv5}{rgb}{0.0,0.0,0.0}
\definecolor{cv6}{rgb}{0.0,0.0,0.0}
\definecolor{cfv6}{rgb}{1.0,1.0,1.0}
\definecolor{clv6}{rgb}{0.0,0.0,0.0}
\definecolor{cv0v1}{rgb}{0.0,0.0,0.0}
\definecolor{cv0v2}{rgb}{0.0,0.0,0.0}
\definecolor{cv0v5}{rgb}{0.0,0.0,0.0}
\definecolor{cv0v6}{rgb}{0.0,0.0,0.0}
\definecolor{cv1v2}{rgb}{0.0,0.0,0.0}
\definecolor{cv1v3}{rgb}{0.0,0.0,0.0}
\definecolor{cv1v6}{rgb}{0.0,0.0,0.0}
\definecolor{cv2v3}{rgb}{0.0,0.0,0.0}
\definecolor{cv2v4}{rgb}{0.0,0.0,0.0}
\definecolor{cv3v4}{rgb}{0.0,0.0,0.0}
\definecolor{cv3v5}{rgb}{0.0,0.0,0.0}
\definecolor{cv4v5}{rgb}{0.0,0.0,0.0}
\definecolor{cv4v6}{rgb}{0.0,0.0,0.0}
\definecolor{cv5v6}{rgb}{0.0,0.0,0.0}
\Vertex[style={minimum size=1.0cm,draw=cv0,fill=cfv0,text=clv0,shape=circle},LabelOut=false,L=\hbox{$0$},x=2.5cm,y=5.0cm]{v0}
\Vertex[style={minimum size=1.0cm,draw=cv1,fill=cfv1,text=clv1,shape=circle},LabelOut=false,L=\hbox{$1$},x=0.4952cm,y=4.0097cm]{v1}
\Vertex[style={minimum size=1.0cm,draw=cv2,fill=cfv2,text=clv2,shape=circle},LabelOut=false,L=\hbox{$2$},x=0.0cm,y=1.7845cm]{v2}
\Vertex[style={minimum size=1.0cm,draw=cv3,fill=cfv3,text=clv3,shape=circle},LabelOut=false,L=\hbox{$3$},x=1.3874cm,y=0.0cm]{v3}
\Vertex[style={minimum size=1.0cm,draw=cv4,fill=cfv4,text=clv4,shape=circle},LabelOut=false,L=\hbox{$4$},x=3.6126cm,y=0.0cm]{v4}
\Vertex[style={minimum size=1.0cm,draw=cv5,fill=cfv5,text=clv5,shape=circle},LabelOut=false,L=\hbox{$5$},x=5.0cm,y=1.7845cm]{v5}
\Vertex[style={minimum size=1.0cm,draw=cv6,fill=cfv6,text=clv6,shape=circle},LabelOut=false,L=\hbox{$6$},x=4.5048cm,y=4.0097cm]{v6}
\Edge[lw=0.1cm,style={color=cv0v1,},](v0)(v1)
\Edge[lw=0.1cm,style={color=cv0v2,},](v0)(v2)
\Edge[lw=0.1cm,style={color=cv0v5,},](v0)(v5)
\Edge[lw=0.1cm,style={color=cv0v6,},](v0)(v6)
\Edge[lw=0.1cm,style={color=cv1v2,},](v1)(v2)
\Edge[lw=0.1cm,style={color=cv1v3,},](v1)(v3)
\Edge[lw=0.1cm,style={color=cv1v6,},](v1)(v6)
\Edge[lw=0.1cm,style={color=cv2v3,},](v2)(v3)
\Edge[lw=0.1cm,style={color=cv2v4,},](v2)(v4)
\Edge[lw=0.1cm,style={color=cv3v4,},](v3)(v4)
\Edge[lw=0.1cm,style={color=cv3v5,},](v3)(v5)
\Edge[lw=0.1cm,style={color=cv4v5,},](v4)(v5)
\Edge[lw=0.1cm,style={color=cv4v6,},](v4)(v6)
\Edge[lw=0.1cm,style={color=cv5v6,},](v5)(v6)
\end{tikzpicture}
\caption{The $C_7^2$ graph}
\end{figure}
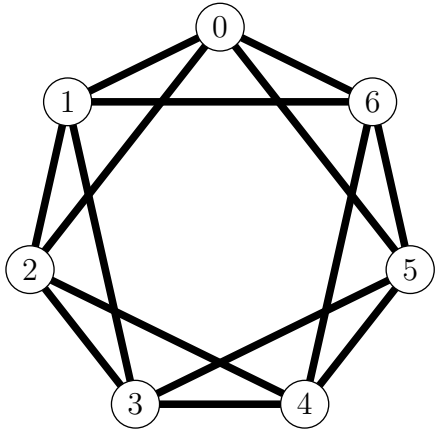

Its sandpile group is non-cyclic. Our goal is to show that with the idea mentioned above it can be solved efficiently. Furthermore using the special structure of the group, one can give parametric solution for the problem, making it even more vulnerable. 

In \cite{cn2_sandpile} the sandpile group of $C_n^2$, denoted by $\mathcal{S}(C_n^2)$ is described. It is the direct sum of two or three cyclic groups depending on the parameter  $n$, more precisely 
$$\mathcal{S}(C_n^2)\cong \Z_{(n,F_n)} \oplus \Z_{F_n} \oplus \Z_{\frac{nF_n}{(n,F_n)}},$$ where $F_n$ denotes the $n$th Fibonacci number. Considering the Laplacian matrix of the $C_n^2$ graph it can be easily seen that it is a circulant matrix. Cline, Plemmons and Worm showed in \cite{generalized_inverse} that the pseudoinverse of these type of matrices are also circulant. Moreover they could give a computational form for the pseudoinverse. Combining these results we can give parametric solution to the DLP based on the input factors. 

Let $c_1, c_2 \in \mathcal{S}(C_n^2)$ be configurations, which are the elements of the sandpile group. The DLP can be described as the following: we are looking for a $2 \leq x \leq \mbox{ord}(\mathcal{S}(C_n^2))$ such that $(x \cdot c_1)^{\circ} = c_2$, where $(x \cdot c_1)^{\circ}$ means the stabilization of the configuration in the group.
 Denote the Laplacian of $\mathcal{S}(C_n^2)$ with $$L = \begin{pmatrix}
        a_0 & a_1 & \hdots & a_{n-1} \\
        a_{n-1} & a_0  & \hdots & a_{n-2} \\
        \vdots & \vdots  & \ddots & \vdots \\
        a_{1} & a_{2}  & \hdots & a_{0} \\
        \end{pmatrix}.$$

The main steps of solving the DLP are the following:
\begin{itemize}
    \item We can calculate the pseudoinverse, using the form given by Cline, Plemmons and Worm. Let $\mu_0, \hdots,\mu_{n-1}$ the eigenvalues, $\omega$ primitive $n$th root of unity, $\lambda$ is also an $n$th root of $1$ in this special case, because the Laplacian $L$ is a $k$-circulant matrix with $k=1.$
Keeping the notation from $L,$ the  first row of the Laplacian's pseudoinverse matrix $(b_0, \hdots, b_{n-1})$ can be given by: 
    $$b_i = \frac{1}{n}\sum_{j = 0}^{n-1} \beta_j ( \lambda \omega^j)^{-i}, \quad i = 0,1, \hdots, n-1,$$
    where
    \begin{equation*}
    \beta_j=
    \begin{cases}
    0 & \mbox{ if }\mu_j=0,\\
    \frac{1}{\mu_j} & \mbox{ if }\mu_j\neq 0.
    \end{cases}
    \end{equation*}
 \item Denote the pseudoinverse with $P$, the divisors of the configurations $c_i$ with $\overline{c_i}$, the generators of the group with $g_j$, their divisors with $\overline{g_j}$, calculate $\overline{c_1}^T \cdot P \cdot \overline{g_j} = r_{j,1} + \Z$, and $\overline{c_2}^T \cdot P \cdot \overline{g_j} = r_{j, 2} + \Z$.
  \item Solve the Diophantine equations $r_{j,2} = r_{j,1} x + y.$ Using the Chinese remainder theorem we get $x$ modulo the least common multiple of the denominators of the rational numbers $r_{j,1}$. 
\end{itemize}

Let us take a simple example to see details of the computations we need to perform. Let $n=7.$ The graph is the  one given by Figure 1. The divisors corresponding to the generators are as follows
\begin{eqnarray*}
\overline{g_1}&=&(-10, 1, 1, 3, 2, 3, 0),\\
\overline{g_2}&=&(-13, 3, 1, 3, 3, 1, 2).
\end{eqnarray*}
The Laplacian matrix is
$$
\left(\begin{array}{rrrrrrr}
4 & -1 & -1 & 0 & 0 & -1 & -1 \\
-1 & 4 & -1 & -1 & 0 & 0 & -1 \\
-1 & -1 & 4 & -1 & -1 & 0 & 0 \\
0 & -1 & -1 & 4 & -1 & -1 & 0 \\
0 & 0 & -1 & -1 & 4 & -1 & -1 \\
-1 & 0 & 0 & -1 & -1 & 4 & -1 \\
-1 & -1 & 0 & 0 & -1 & -1 & 4
\end{array}\right)
$$
and the pseudoinverse matrix is given by
$$
P=
\left(\begin{array}{rrrrrrr}
\frac{18}{91} & -\frac{1}{91} & -\frac{2}{91} & -\frac{6}{91} & -\frac{6}{91} & -\frac{2}{91} & -\frac{1}{91} \\
-\frac{1}{91} & \frac{18}{91} & -\frac{1}{91} & -\frac{2}{91} & -\frac{6}{91} & -\frac{6}{91} & -\frac{2}{91} \\
-\frac{2}{91} & -\frac{1}{91} & \frac{18}{91} & -\frac{1}{91} & -\frac{2}{91} & -\frac{6}{91} & -\frac{6}{91} \\
-\frac{6}{91} & -\frac{2}{91} & -\frac{1}{91} & \frac{18}{91} & -\frac{1}{91} & -\frac{2}{91} & -\frac{6}{91} \\
-\frac{6}{91} & -\frac{6}{91} & -\frac{2}{91} & -\frac{1}{91} & \frac{18}{91} & -\frac{1}{91} & -\frac{2}{91} \\
-\frac{2}{91} & -\frac{6}{91} & -\frac{6}{91} & -\frac{2}{91} & -\frac{1}{91} & \frac{18}{91} & -\frac{1}{91} \\
-\frac{1}{91} & -\frac{2}{91} & -\frac{6}{91} & -\frac{6}{91} & -\frac{2}{91} & -\frac{1}{91} & \frac{18}{91}
\end{array}\right).
$$
Consider the DLP $xc_1=c_2,$ where the divisors corresponding to $c_1$ and $c_2$ are
\begin{eqnarray*}
\overline{c_1}&=&(-12, 3, 3, 0, 3, 3, 0),\\
\overline{c_2}&=&(-14, 3, 3, 3, 3, 1, 1).
\end{eqnarray*}
We obtain that 
\begin{eqnarray*}
\overline{c_1}^T\cdot P\cdot \overline{g_1}&=&435/13 \longrightarrow r_{1,1}=6/13,\\
\overline{c_1}^T\cdot P\cdot \overline{g_2}&=&3825/91 \longrightarrow r_{2,1}=3/91,\\
\overline{c_2}^T\cdot P\cdot \overline{g_1}&=&523/13 \longrightarrow r_{1,2}=3/13,\\
\overline{c_2}^T\cdot P\cdot \overline{g_2}&=&4701/91 \longrightarrow r_{2,2}=60/91.
\end{eqnarray*}
The solutions of the linear Diophantine equation $3=6x+13y$ are given by
$$
x=-6+13t,\qquad y=3-6t,
$$
hence $x\equiv 7\pmod{13}.$
The second linear Diophantine equation is $60=3x+91y$ having the solutions
$$
x=-1800+91t,\qquad y=60-3t,
$$
therefore $x\equiv 20\pmod{91}.$ Applying the Chinese remainder theorem yields that $x\equiv 20\pmod{91}.$

\section{DLP in wheel graphs}
The wheel graph $W_n$ has $n+1$ vertices, the graph contains a cycle of length $n$ and there is an additional vertex adjacent to all vertices of the cycle.
\begin{figure}[h]
\begin{tikzpicture}
\definecolor{cv0}{rgb}{0.0,0.0,0.0}
\definecolor{cfv0}{rgb}{1.0,1.0,1.0}
\definecolor{clv0}{rgb}{0.0,0.0,0.0}
\definecolor{cv1}{rgb}{0.0,0.0,0.0}
\definecolor{cfv1}{rgb}{1.0,1.0,1.0}
\definecolor{clv1}{rgb}{0.0,0.0,0.0}
\definecolor{cv2}{rgb}{0.0,0.0,0.0}
\definecolor{cfv2}{rgb}{1.0,1.0,1.0}
\definecolor{clv2}{rgb}{0.0,0.0,0.0}
\definecolor{cv3}{rgb}{0.0,0.0,0.0}
\definecolor{cfv3}{rgb}{1.0,1.0,1.0}
\definecolor{clv3}{rgb}{0.0,0.0,0.0}
\definecolor{cv4}{rgb}{0.0,0.0,0.0}
\definecolor{cfv4}{rgb}{1.0,1.0,1.0}
\definecolor{clv4}{rgb}{0.0,0.0,0.0}
\definecolor{cv5}{rgb}{0.0,0.0,0.0}
\definecolor{cfv5}{rgb}{1.0,1.0,1.0}
\definecolor{clv5}{rgb}{0.0,0.0,0.0}
\definecolor{cv6}{rgb}{0.0,0.0,0.0}
\definecolor{cfv6}{rgb}{1.0,1.0,1.0}
\definecolor{clv6}{rgb}{0.0,0.0,0.0}
\definecolor{cv7}{rgb}{0.0,0.0,0.0}
\definecolor{cfv7}{rgb}{1.0,1.0,1.0}
\definecolor{clv7}{rgb}{0.0,0.0,0.0}
\definecolor{cv0v1}{rgb}{0.0,0.0,0.0}
\definecolor{cv0v6}{rgb}{0.0,0.0,0.0}
\definecolor{cv7v0}{rgb}{0.0,0.0,0.0}
\definecolor{cv1v2}{rgb}{0.0,0.0,0.0}
\definecolor{cv7v1}{rgb}{0.0,0.0,0.0}
\definecolor{cv2v3}{rgb}{0.0,0.0,0.0}
\definecolor{cv7v2}{rgb}{0.0,0.0,0.0}
\definecolor{cv3v4}{rgb}{0.0,0.0,0.0}
\definecolor{cv7v3}{rgb}{0.0,0.0,0.0}
\definecolor{cv4v5}{rgb}{0.0,0.0,0.0}
\definecolor{cv7v4}{rgb}{0.0,0.0,0.0}
\definecolor{cv5v6}{rgb}{0.0,0.0,0.0}
\definecolor{cv7v5}{rgb}{0.0,0.0,0.0}
\definecolor{cv7v6}{rgb}{0.0,0.0,0.0}
\Vertex[style={minimum size=1.0cm,draw=cv0,fill=cfv0,text=clv0,shape=circle},LabelOut=false,L=\hbox{$1$},x=2.5cm,y=5.0cm]{v0}
\Vertex[style={minimum size=1.0cm,draw=cv1,fill=cfv1,text=clv1,shape=circle},LabelOut=false,L=\hbox{$2$},x=0.4952cm,y=4.0097cm]{v1}
\Vertex[style={minimum size=1.0cm,draw=cv2,fill=cfv2,text=clv2,shape=circle},LabelOut=false,L=\hbox{$3$},x=0.0cm,y=1.7845cm]{v2}
\Vertex[style={minimum size=1.0cm,draw=cv3,fill=cfv3,text=clv3,shape=circle},LabelOut=false,L=\hbox{$4$},x=1.3874cm,y=0.0cm]{v3}
\Vertex[style={minimum size=1.0cm,draw=cv4,fill=cfv4,text=clv4,shape=circle},LabelOut=false,L=\hbox{$5$},x=3.6126cm,y=0.0cm]{v4}
\Vertex[style={minimum size=1.0cm,draw=cv5,fill=cfv5,text=clv5,shape=circle},LabelOut=false,L=\hbox{$6$},x=5.0cm,y=1.7845cm]{v5}
\Vertex[style={minimum size=1.0cm,draw=cv6,fill=cfv6,text=clv6,shape=circle},LabelOut=false,L=\hbox{$7$},x=4.5048cm,y=4.0097cm]{v6}
\Vertex[style={minimum size=1.0cm,draw=cv7,fill=cfv7,text=clv7,shape=circle},LabelOut=false,L=\hbox{$0$},x=2.5cm,y=2.3698cm]{v7}
\Edge[lw=0.1cm,style={color=cv0v1,},](v0)(v1)
\Edge[lw=0.1cm,style={color=cv0v6,},](v0)(v6)
\Edge[lw=0.1cm,style={color=cv7v0,},](v7)(v0)
\Edge[lw=0.1cm,style={color=cv1v2,},](v1)(v2)
\Edge[lw=0.1cm,style={color=cv7v1,},](v7)(v1)
\Edge[lw=0.1cm,style={color=cv2v3,},](v2)(v3)
\Edge[lw=0.1cm,style={color=cv7v2,},](v7)(v2)
\Edge[lw=0.1cm,style={color=cv3v4,},](v3)(v4)
\Edge[lw=0.1cm,style={color=cv7v3,},](v7)(v3)
\Edge[lw=0.1cm,style={color=cv4v5,},](v4)(v5)
\Edge[lw=0.1cm,style={color=cv7v4,},](v7)(v4)
\Edge[lw=0.1cm,style={color=cv5v6,},](v5)(v6)
\Edge[lw=0.1cm,style={color=cv7v5,},](v7)(v5)
\Edge[lw=0.1cm,style={color=cv7v6,},](v7)(v6)
\end{tikzpicture}
\caption{The $W_7$ graph}
\end{figure}
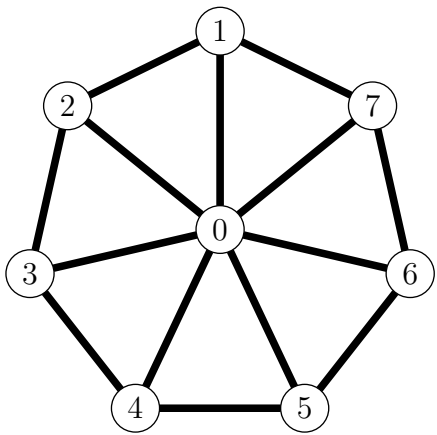
\begin{center}
\begin{tabular}{|r|rrrrrrrr|}
\hline
$n$ & $3$ & $4$ & $5$ & $6$ & $7$ & $8$ & $9$ & $10$ \\
tree-number & $16$ & $45$ & $121$ & $320$ & $841$ & $2205$ & $5776$ & $15125$ \\
\hline
\end{tabular}
\end{center}
One can easily see the two patterns of the tree-numbers in the above table, Biggs \cite{GWheel} proved that in case of odd values of $n$ the tree-number is a square of a Lucas number and in case of even values of $n$ it is 5 times a square of a Fibonacci number. If $n$ is odd, let say $n=2r+1,$ then Biggs also provided generators of the sandpile group. One denotes the vertices of the rim of the wheel by 
$$-r,-r+1,\ldots,-1,0,1,\ldots,r-1,r\pmod{n}$$ and define a configuration $g$ as follows
\begin{equation*}
g(v)=
\begin{cases}
1, & \mbox{ if } v=\pm r;\\
2, & \mbox{ otherwise.}
\end{cases}
\end{equation*}
Two generators are given by $g_1=g$ and $g_2=g(v-1).$ Once the pseudoinverse is computed of the Laplacian matrix, the DLP can be computed efficiently. We provide a simple example to see further details. We consider the case $n=7,$ so we deal with the graph $W_7.$ The Laplacian matrix is given by 
$$
\left(\begin{array}{rrrrrrrr}
7 & -1 & -1 & -1 & -1 & -1 & -1 & -1 \\
-1 & 3 & -1 & 0 & 0 & 0 & 0 & -1 \\
-1 & -1 & 3 & -1 & 0 & 0 & 0 & 0 \\
-1 & 0 & -1 & 3 & -1 & 0 & 0 & 0 \\
-1 & 0 & 0 & -1 & 3 & -1 & 0 & 0 \\
-1 & 0 & 0 & 0 & -1 & 3 & -1 & 0 \\
-1 & 0 & 0 & 0 & 0 & -1 & 3 & -1 \\
-1 & -1 & 0 & 0 & 0 & 0 & -1 & 3
\end{array}\right)
$$ 
and the pseudoinverse is as follows
$$
P=
\left(\begin{array}{rrrrrrrr}
\frac{7}{64} & -\frac{1}{64} & -\frac{1}{64} & -\frac{1}{64} & -\frac{1}{64} & -\frac{1}{64} & -\frac{1}{64} & -\frac{1}{64} \\
-\frac{1}{64} & \frac{571}{1856} & \frac{59}{1856} & -\frac{133}{1856} & -\frac{197}{1856} & -\frac{197}{1856} & -\frac{133}{1856} & \frac{59}{1856} \\
-\frac{1}{64} & \frac{59}{1856} & \frac{571}{1856} & \frac{59}{1856} & -\frac{133}{1856} & -\frac{197}{1856} & -\frac{197}{1856} & -\frac{133}{1856} \\
-\frac{1}{64} & -\frac{133}{1856} & \frac{59}{1856} & \frac{571}{1856} & \frac{59}{1856} & -\frac{133}{1856} & -\frac{197}{1856} & -\frac{197}{1856} \\
-\frac{1}{64} & -\frac{197}{1856} & -\frac{133}{1856} & \frac{59}{1856} & \frac{571}{1856} & \frac{59}{1856} & -\frac{133}{1856} & -\frac{197}{1856} \\
-\frac{1}{64} & -\frac{197}{1856} & -\frac{197}{1856} & -\frac{133}{1856} & \frac{59}{1856} & \frac{571}{1856} & \frac{59}{1856} & -\frac{133}{1856} \\
-\frac{1}{64} & -\frac{133}{1856} & -\frac{197}{1856} & -\frac{197}{1856} & -\frac{133}{1856} & \frac{59}{1856} & \frac{571}{1856} & \frac{59}{1856} \\
-\frac{1}{64} & \frac{59}{1856} & -\frac{133}{1856} & -\frac{197}{1856} & -\frac{197}{1856} & -\frac{133}{1856} & \frac{59}{1856} & \frac{571}{1856}
\end{array}\right).
$$
The sandpile group of the graph is $\mathcal{S}(W_7)\cong \Z_{L_7}\oplus \Z_{L_7}$ and the generators provided by Biggs are as follows
$$
g_1=(1,2,2,2,2,2,1),\quad g_2=(1,1,2,2,2,2,2).
$$
Consider the DLP $xc_1=c_2,$ where $c_1=(2, 2, 2, 0, 2, 2, 0)$ and $c_2=(2, 2, 2, 1, 2, 2, 1).$
We use the following divisors to compute $r_{j,1}$ and $r_{j,2},j=1,2:$
\begin{eqnarray*}
\overline{g_1}=(-12,1,2,2,2,2,2,1), && \overline{g_2}=(-12,1,1,2,2,2,2,2),\\
\overline{c_1}=(-10,2, 2, 2, 0, 2, 2, 0), && \overline{c_2}=(-12,2, 2, 2, 1, 2, 2, 1).
\end{eqnarray*} 
We get that
\begin{eqnarray*}
\overline{c_1}^T\cdot P\cdot \overline{g_1}&=&504/29 \longrightarrow r_{1,1}=11/29,\\
\overline{c_1}^T\cdot P\cdot \overline{g_2}&=&484/29 \longrightarrow r_{2,1}=20/29,\\
\overline{c_2}^T\cdot P\cdot \overline{g_1}&=&600/29 \longrightarrow r_{1,2}=20/29,\\
\overline{c_2}^T\cdot P\cdot \overline{g_2}&=&590/29 \longrightarrow r_{2,2}=10/29.
\end{eqnarray*}
So we need to solve the linear Diophantine equations
$$
20=11x+29y,\qquad 10=20x+29y.
$$
The solution of the first one has $x\equiv 15\pmod{29}$ and it turns out that the second has the same solution. Therefore $x\equiv 15\pmod{29}.$

We implemented the procedure in SageMath \cite{sage} to check the CPU times in case of larger groups, in the table we indicate $n,$ the order of the group and the CPU time (in millisecond) of a randomly generated DLP for given $n$
\begin{center}
\begin{tabular}{|r|r|r|}
\hline
$n$ & order of the group & CPU time\\
\hline
29 & 1322157322201 & 51ms\\
31 & 9062201101801 & 55ms\\
33 & 62113250390416 & 58ms\\
35 & 425730551631121 & 67ms\\
37 & 2918000611027441 & 69ms\\
39 & 20000273725560976 & 72ms\\
41 & 137083915467899401 & 79ms\\
43 & 939587134549734841 & 82ms\\
45 & 6440026026380244496 & 89ms\\
\hline
\end{tabular}
\end{center}

\section{DLP in subdivided banana graphs}
Consider a tuple of positive integers $s=(s_1,s_2,\ldots,s_m).$ A banana graph $B_m$ has two vertices and $m$ edges between them. The $s$-subdivided banana graph is obtained from $B_m$ by subdividing the $i$th edge $s_i -1$ times. 
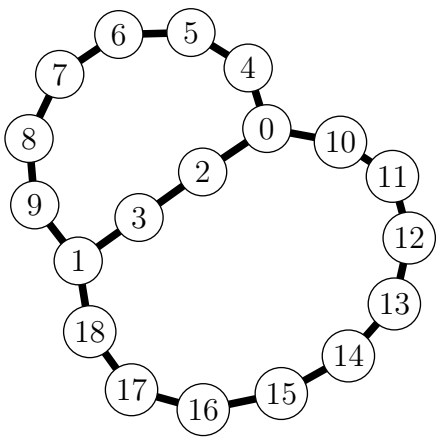
\begin{figure}[h]
\begin{tikzpicture}
\definecolor{cv0}{rgb}{0.0,0.0,0.0}
\definecolor{cfv0}{rgb}{1.0,1.0,1.0}
\definecolor{clv0}{rgb}{0.0,0.0,0.0}
\definecolor{cv1}{rgb}{0.0,0.0,0.0}
\definecolor{cfv1}{rgb}{1.0,1.0,1.0}
\definecolor{clv1}{rgb}{0.0,0.0,0.0}
\definecolor{cv2}{rgb}{0.0,0.0,0.0}
\definecolor{cfv2}{rgb}{1.0,1.0,1.0}
\definecolor{clv2}{rgb}{0.0,0.0,0.0}
\definecolor{cv3}{rgb}{0.0,0.0,0.0}
\definecolor{cfv3}{rgb}{1.0,1.0,1.0}
\definecolor{clv3}{rgb}{0.0,0.0,0.0}
\definecolor{cv4}{rgb}{0.0,0.0,0.0}
\definecolor{cfv4}{rgb}{1.0,1.0,1.0}
\definecolor{clv4}{rgb}{0.0,0.0,0.0}
\definecolor{cv5}{rgb}{0.0,0.0,0.0}
\definecolor{cfv5}{rgb}{1.0,1.0,1.0}
\definecolor{clv5}{rgb}{0.0,0.0,0.0}
\definecolor{cv6}{rgb}{0.0,0.0,0.0}
\definecolor{cfv6}{rgb}{1.0,1.0,1.0}
\definecolor{clv6}{rgb}{0.0,0.0,0.0}
\definecolor{cv7}{rgb}{0.0,0.0,0.0}
\definecolor{cfv7}{rgb}{1.0,1.0,1.0}
\definecolor{clv7}{rgb}{0.0,0.0,0.0}
\definecolor{cv8}{rgb}{0.0,0.0,0.0}
\definecolor{cfv8}{rgb}{1.0,1.0,1.0}
\definecolor{clv8}{rgb}{0.0,0.0,0.0}
\definecolor{cv9}{rgb}{0.0,0.0,0.0}
\definecolor{cfv9}{rgb}{1.0,1.0,1.0}
\definecolor{clv9}{rgb}{0.0,0.0,0.0}
\definecolor{cv10}{rgb}{0.0,0.0,0.0}
\definecolor{cfv10}{rgb}{1.0,1.0,1.0}
\definecolor{clv10}{rgb}{0.0,0.0,0.0}
\definecolor{cv11}{rgb}{0.0,0.0,0.0}
\definecolor{cfv11}{rgb}{1.0,1.0,1.0}
\definecolor{clv11}{rgb}{0.0,0.0,0.0}
\definecolor{cv12}{rgb}{0.0,0.0,0.0}
\definecolor{cfv12}{rgb}{1.0,1.0,1.0}
\definecolor{clv12}{rgb}{0.0,0.0,0.0}
\definecolor{cv13}{rgb}{0.0,0.0,0.0}
\definecolor{cfv13}{rgb}{1.0,1.0,1.0}
\definecolor{clv13}{rgb}{0.0,0.0,0.0}
\definecolor{cv14}{rgb}{0.0,0.0,0.0}
\definecolor{cfv14}{rgb}{1.0,1.0,1.0}
\definecolor{clv14}{rgb}{0.0,0.0,0.0}
\definecolor{cv15}{rgb}{0.0,0.0,0.0}
\definecolor{cfv15}{rgb}{1.0,1.0,1.0}
\definecolor{clv15}{rgb}{0.0,0.0,0.0}
\definecolor{cv16}{rgb}{0.0,0.0,0.0}
\definecolor{cfv16}{rgb}{1.0,1.0,1.0}
\definecolor{clv16}{rgb}{0.0,0.0,0.0}
\definecolor{cv17}{rgb}{0.0,0.0,0.0}
\definecolor{cfv17}{rgb}{1.0,1.0,1.0}
\definecolor{clv17}{rgb}{0.0,0.0,0.0}
\definecolor{cv18}{rgb}{0.0,0.0,0.0}
\definecolor{cfv18}{rgb}{1.0,1.0,1.0}
\definecolor{clv18}{rgb}{0.0,0.0,0.0}
\definecolor{cv0v2}{rgb}{0.0,0.0,0.0}
\definecolor{cv0v4}{rgb}{0.0,0.0,0.0}
\definecolor{cv0v10}{rgb}{0.0,0.0,0.0}
\definecolor{cv1v18}{rgb}{0.0,0.0,0.0}
\definecolor{cv1v3}{rgb}{0.0,0.0,0.0}
\definecolor{cv1v9}{rgb}{0.0,0.0,0.0}
\definecolor{cv2v3}{rgb}{0.0,0.0,0.0}
\definecolor{cv4v5}{rgb}{0.0,0.0,0.0}
\definecolor{cv5v6}{rgb}{0.0,0.0,0.0}
\definecolor{cv6v7}{rgb}{0.0,0.0,0.0}
\definecolor{cv7v8}{rgb}{0.0,0.0,0.0}
\definecolor{cv8v9}{rgb}{0.0,0.0,0.0}
\definecolor{cv10v11}{rgb}{0.0,0.0,0.0}
\definecolor{cv11v12}{rgb}{0.0,0.0,0.0}
\definecolor{cv12v13}{rgb}{0.0,0.0,0.0}
\definecolor{cv13v14}{rgb}{0.0,0.0,0.0}
\definecolor{cv14v15}{rgb}{0.0,0.0,0.0}
\definecolor{cv15v16}{rgb}{0.0,0.0,0.0}
\definecolor{cv16v17}{rgb}{0.0,0.0,0.0}
\definecolor{cv17v18}{rgb}{0.0,0.0,0.0}
\Vertex[style={minimum size=1.0cm,draw=cv0,fill=cfv0,text=clv0,shape=circle},LabelOut=false,L=\hbox{$0$},x=3.1274cm,y=3.7268cm]{v0}
\Vertex[style={minimum size=1.0cm,draw=cv1,fill=cfv1,text=clv1,shape=circle},LabelOut=false,L=\hbox{$1$},x=0.6465cm,y=1.9716cm]{v1}
\Vertex[style={minimum size=1.0cm,draw=cv2,fill=cfv2,text=clv2,shape=circle},LabelOut=false,L=\hbox{$2$},x=2.2833cm,y=3.1514cm]{v2}
\Vertex[style={minimum size=1.0cm,draw=cv3,fill=cfv3,text=clv3,shape=circle},LabelOut=false,L=\hbox{$3$},x=1.4468cm,y=2.5475cm]{v3}
\Vertex[style={minimum size=1.0cm,draw=cv4,fill=cfv4,text=clv4,shape=circle},LabelOut=false,L=\hbox{$4$},x=2.882cm,y=4.5307cm]{v4}
\Vertex[style={minimum size=1.0cm,draw=cv5,fill=cfv5,text=clv5,shape=circle},LabelOut=false,L=\hbox{$5$},x=2.1291cm,y=5.0cm]{v5}
\Vertex[style={minimum size=1.0cm,draw=cv6,fill=cfv6,text=clv6,shape=circle},LabelOut=false,L=\hbox{$6$},x=1.181cm,y=4.9742cm]{v6}
\Vertex[style={minimum size=1.0cm,draw=cv7,fill=cfv7,text=clv7,shape=circle},LabelOut=false,L=\hbox{$7$},x=0.4034cm,y=4.4459cm]{v7}
\Vertex[style={minimum size=1.0cm,draw=cv8,fill=cfv8,text=clv8,shape=circle},LabelOut=false,L=\hbox{$8$},x=0.0cm,y=3.5946cm]{v8}
\Vertex[style={minimum size=1.0cm,draw=cv9,fill=cfv9,text=clv9,shape=circle},LabelOut=false,L=\hbox{$9$},x=0.0749cm,y=2.714cm]{v9}
\Vertex[style={minimum size=1.0cm,draw=cv10,fill=cfv10,text=clv10,shape=circle},LabelOut=false,L=\hbox{$10$},x=4.0957cm,y=3.547cm]{v10}
\Vertex[style={minimum size=1.0cm,draw=cv11,fill=cfv11,text=clv11,shape=circle},LabelOut=false,L=\hbox{$11$},x=4.7792cm,y=3.0932cm]{v11}
\Vertex[style={minimum size=1.0cm,draw=cv12,fill=cfv12,text=clv12,shape=circle},LabelOut=false,L=\hbox{$12$},x=5.0cm,y=2.2757cm]{v12}
\Vertex[style={minimum size=1.0cm,draw=cv13,fill=cfv13,text=clv13,shape=circle},LabelOut=false,L=\hbox{$13$},x=4.8046cm,y=1.4015cm]{v13}
\Vertex[style={minimum size=1.0cm,draw=cv14,fill=cfv14,text=clv14,shape=circle},LabelOut=false,L=\hbox{$14$},x=4.2001cm,y=0.7117cm]{v14}
\Vertex[style={minimum size=1.0cm,draw=cv15,fill=cfv15,text=clv15,shape=circle},LabelOut=false,L=\hbox{$15$},x=3.3176cm,y=0.2137cm]{v15}
\Vertex[style={minimum size=1.0cm,draw=cv16,fill=cfv16,text=clv16,shape=circle},LabelOut=false,L=\hbox{$16$},x=2.2912cm,y=0.0cm]{v16}
\Vertex[style={minimum size=1.0cm,draw=cv17,fill=cfv17,text=clv17,shape=circle},LabelOut=false,L=\hbox{$17$},x=1.3479cm,y=0.2634cm]{v17}
\Vertex[style={minimum size=1.0cm,draw=cv18,fill=cfv18,text=clv18,shape=circle},LabelOut=false,L=\hbox{$18$},x=0.7984cm,y=1.0335cm]{v18}
\Edge[lw=0.1cm,style={color=cv0v2,},](v0)(v2)
\Edge[lw=0.1cm,style={color=cv0v4,},](v0)(v4)
\Edge[lw=0.1cm,style={color=cv0v10,},](v0)(v10)
\Edge[lw=0.1cm,style={color=cv1v18,},](v1)(v18)
\Edge[lw=0.1cm,style={color=cv1v3,},](v1)(v3)
\Edge[lw=0.1cm,style={color=cv1v9,},](v1)(v9)
\Edge[lw=0.1cm,style={color=cv2v3,},](v2)(v3)
\Edge[lw=0.1cm,style={color=cv4v5,},](v4)(v5)
\Edge[lw=0.1cm,style={color=cv5v6,},](v5)(v6)
\Edge[lw=0.1cm,style={color=cv6v7,},](v6)(v7)
\Edge[lw=0.1cm,style={color=cv7v8,},](v7)(v8)
\Edge[lw=0.1cm,style={color=cv8v9,},](v8)(v9)
\Edge[lw=0.1cm,style={color=cv10v11,},](v10)(v11)
\Edge[lw=0.1cm,style={color=cv11v12,},](v11)(v12)
\Edge[lw=0.1cm,style={color=cv12v13,},](v12)(v13)
\Edge[lw=0.1cm,style={color=cv13v14,},](v13)(v14)
\Edge[lw=0.1cm,style={color=cv14v15,},](v14)(v15)
\Edge[lw=0.1cm,style={color=cv15v16,},](v15)(v16)
\Edge[lw=0.1cm,style={color=cv16v17,},](v16)(v17)
\Edge[lw=0.1cm,style={color=cv17v18,},](v17)(v18)
\end{tikzpicture}
\caption{The $B_{(3,7,10)}$ graph}
\end{figure}
In some cases the sandpile group of subdivided banana graphs can be described and it turns out that the group is cyclic, the following result is from \cite{Banana}.
\begin{proposition}
Fix a prime $p$ and an integer $r.$ Let $s = (s_1,s_2,\ldots,s_m )$ be a tuple of positive
integers such that 
$$
\sum_{i=1}^m \frac{\prod_{j=1}^m s_j}{s_i}=p^r
$$
and $\gcd(s_i , p) = 1$ for all $i.$ Then
$$\mbox{Jac}(B_s)\cong (\mathbb{Z}/p^r\mathbb{Z},\langle \cdot,\cdot\rangle),$$
where  $\langle \cdot,\cdot\rangle$ is the pairing on $\mathbb{Z}/p^r\mathbb{Z}$ given by
$$
\langle x,y\rangle=\frac{\left(\prod_{i=1}^m s_i\right) xy}{p^r}.
$$
\end{proposition}
In the proof a generator of the sandpile group is also provided in these cases. 
If we denote the two vertices of the banana graph $B_m$ by $v_0$ and $v_1,$ then the divisor $D=v_0-v_1$ generates the sandpile group of $B_s.$ In these cases we not only know generators in a parametric form, but also the pairing helps to obtain the linear Diophantine equation without computing the pseudoinverse of the Laplacian. Let us see a simple example given at Figure 3, that is the subdivided banana graph $B_{(3,7,10)}.$ 
Consider the DLP with
\begin{eqnarray*}
\overline{c_1}&=&(2, -17, 1, 1, 1, 1, 1, 1, 1, 0, 1, 1, 0, 1, 1, 1, 1, 1, 1),\\
\overline{c_2}&=&(2, -17, 1, 1, 1, 1, 0, 1, 1, 1, 1, 1, 1, 1, 1, 1, 0, 1, 1).
\end{eqnarray*}
The values of the pairing can be determined by the function $f: V(G)\rightarrow \mathbb{Z}$ for which 
$$b(e)=b_i=\frac{\prod_{j=1}^m s_j}{s_i}=f(e_h)-f(e_t)$$
for any edge $e$ in the $i$th subdivided edge, where $e_h$ denotes the head, $e_t$ denotes the tail of the edge (each edge points toward $v_1$). The order of the group is 121, hence we have
$$
\langle D,D_1 \rangle = \frac{1}{121}\sum_{v\in V(G)}D_1(v)f(v).
$$
We obtain that 
$$
\langle D,c_1 \rangle=\frac{95}{121}, \qquad \langle D,c_2 \rangle=\frac{62}{121}.
$$
It remains to solve the linear Diophantine equation
$$
62=95x+121y.
$$
We obtain that $x=-868+121t, y=682-95t.$ Thus we have $x\equiv -868\equiv 100\pmod{121}.$
\begin{acknowledgement}
Research supported in part by the NKFIH grants 115479, 128088 and 130909 (Sz.~T.).
\end{acknowledgement}

\bibliography{thesis}

\par\noindent\rule{\textwidth}{0.5pt}
\textit{Keywords:} discrete logarithm problem, sandpile group,  critical group\\
\noindent\rule{\textwidth}{0.5pt}

\end{document}